\documentclass{elsart}
\usepackage{amscd,amssymb,graphics}
\usepackage{mathrsfs} 
\usepackage{amsmath}
\usepackage{txfonts,stmaryrd}
\newtheorem{theorem}{Theorem}[section]
\newtheorem{corollary}[theorem]{Corollary} 

\newtheorem{proposition}[theorem]{Proposition}
\theoremstyle{definition}

\newtheorem{definition}[theorem]{Definition}
\theoremstyle{remark}
\newtheorem{remark}[theorem]{Remark}

\newtheorem{example}[theorem]{Example}

\newcommand{\abs}[1]{\vert#1\vert}
\def\norm#1{\left\Vert#1\right\Vert}
\newcommand{\tri}{\hfill$\triangle$}

\def\R {{\mathbb R}}
\def\Q {{\mathbb Q}}

\def\N{{\mathbb N}}
\def\Ur{{\mathbb U}}

\def\e{{\varepsilon}}

\def\Z {{\mathbb Z}}
\def\s{{\mathbb S}}

\def\Iso{{\mathrm{Iso}\,}}

\def\Homeo{{\mathrm{Homeo}\,}}

\def\numero #1{{\vskip 2mm 
\par {\stepcounter{nomer1}(\arabic{nomer1}) #1}}}
\newcounter{nomer1}

\begin{document}
\begin{frontmatter}

\title{A theorem of Hrushovski--Solecki--Vershik applied to
uniform and coarse embeddings of the Urysohn metric space}

\author{Vladimir G. Pestov}
\address{Department of Mathematics and Statistics, 
University of Ottawa, 585 King Edward Ave., Ottawa, Ontario, Canada K1N 6N5
}

\ead{vpest283@uottawa.ca}

\begin{abstract} 
A theorem proved by Hrushovski for graphs and extended by Solecki and Vershik (independently from each other) to metric spaces leads to a stronger version of ultrahomogeneity of the infinite random graph $R$, the universal Urysohn metric space $\Ur$, and other related objects. We show how the result can be used to
average out uniform and coarse embeddings of $\Ur$ (and its various counterparts) into normed spaces. Sometimes this leads to new embeddings of the same kind that are metric transforms and besides extend to affine representations of various isometry groups. 
As an application of this technique, we show that $\Ur$ admits neither a uniform nor a coarse embedding into a uniformly convex Banach space.
\end{abstract}
\begin{keyword}
Urysohn metric space \sep
coarse embedding \sep
uniform embedding \sep
Hrushovsky--Solecki--Vershik theorem \sep
uniform convexity \sep
group of isometries 

\MSC 22A05\sep 43A07 \sep  51F99 \sep 54E35
\end{keyword}
\end{frontmatter}

\section{Introduction}

A theorem by Hrushovski  \cite{hrushovski} states that every finite graph $\Gamma$ can be embedded (as an induced subgraph) into a finite graph $\Gamma^\prime$ so that each partial isomorphism of $\Gamma$ is a restriction of a global automorphism of $\Gamma^\prime$. Solecki \cite{solecki} and (independently) Vershik \cite{Ver05b} have obtained an analogue of the result for metric spaces: every finite metric space $X$ is contained in a finite metric space $Y$ in such a way that partial isometries of $X$ become restrictions of global self-isometries of $Y$. Solecki has deduced the result from a powerful general theorem of Herwig and Lascar \cite{HL}, while Vershik gave a direct construction. 

Both theorems are particular cases of a statement where $X$ and $Y$ are drawn from a class of metric spaces whose distance values belong to a given convex subset $S$ of some additive subsemigroup of real numbers. The Hrushovski theorem is recovered for the set of values $S=\{0,1,2\}$, while the Solecki--Vershik theorem corresponds to the entire real line. 
We outline a proof of the result in Section \ref{s:hsv}; it is also based on the techniques of Herwig and Laskar, but in a different way from Solecki's article.

Let again $S\subseteq\R$ be a subset as above. Denote by $\Ur_S$ a version of the universal Urysohn metric space with distance values in $S$, characterized by the following properties: completeness, ultrahomogeneity (every isometry between finite subspaces is extended to a global self-isometry of the space), and universality for the class of all metric spaces with distance values in $S$. For $S=\{0,1,2\}$ one obtains in this way the infinite random graph $R$ \cite{Rado,cameron}, while $S=\R$ results in the classical Urysohn metric space $\Ur$, and $S=\Z$ leads to the integer-valued Urysohn space $\Ur_\Z$. Now the Hrushovski--Solecki--Vershik theorem can be viewed as a stronger version of ultrahomogeneity for the space $\Ur_S$: every finite subspace $X$ of $\Ur_S$ is contained in a finite subspace $Y$ of $\Ur_S$ in such a way that partial isometries of $X$ are restrictions of global isometries of $\Ur_S$ taking $Y$ to itself. Simple examples (such as the unit sphere in a Hilbert space) show that the property is in general strictly stronger than ultrahomogeneity.

A range of interesting applications of the Hrushovski--Solecki--Vershik property can be found in the original works \cite{hrushovski,solecki,Ver05b}.
Here we apply the result to analysis of (non)existence of uniform and coarse embeddings of the Urysohn metric space into superreflexive Banach spaces.

Let $X$ and $Y$ be two metric spaces, and let $f\colon X\to Y$ be
an embedding of $X$ into $Y$ as a uniform subspace. This property of $f$ is easily seen to be equivalent to the following: for some non-decreasing functions $\rho_1,\rho_2\colon \R_+\to\R_+$ with $0<\rho_1(x)\leq\rho_2(x)$ for $x>0$ and $\rho_2(x)\to 0$ as $x\to 0$, one has for every $x,y\in X$ 
\[\rho_1(d_X(x,y))\leq d_Y(f(x),f(y))\leq \rho_2(d_X(x,y)).\]
Here the double inequality only needs to hold for values of the distance $d_X(x,y)$ in a sufficiently small neighbourhood of zero. 

If, on the contrary, we are interested in the above property holding for values of the distance $d(x,y)$ in a neighbourhood of {\em infinity,} we arrive at the relatively recent notion of a {\em coarse embedding} of $X$ into $Y$. So is called a function $f\colon X\to Y$ (not necessarily continuous) such that, for some non-decreasing, unbounded functions $\rho_1,\rho_2\colon \R_+\to\R_+$ with $\rho_1\leq\rho_2$ one has for every $x,y\in X$ 
\[\rho_1(d_X(x,y))\leq d_Y(f(x),f(y))\leq \rho_2(d_X(x,y)).\]

Coarse embeddings are monomorphisms in the {\em coarse category} of metric spaces, which we will not define here, referring to the book \cite{roe} instead. In the same book, the reader can find a detailed motivation for the concept (the {\em Baum--Connes conjecture}). 

The first example of a metric space admitting no {\em uniform} embedding into the Hilbert space $\ell^2$ was constructed by Per Enflo \cite{enflo}.
Gromov asked in \cite{gromov91} if every separable metric space can be {\em coarsely} embedded into a Hilbert space. The first counter-example was constructed by Dranishnikov, Gong,  Lafforgue, and Yu \cite{dgly}, who have used a suitable modification of Enflo's construction. Thus, the two notions are without doubt linked at some fundamental level, though the full extent of this link is not quite obvious. Notice that initially Gromov even used the same term ``uniform embedding'' to denote what is now known as a coarse embedding. 

Of particular interest in relation to the Baum--Connes conjecture are coarse embeddings of metric
{\em spaces of bounded geometry} (for every $R>0$, the cardinality of each ball $B_R(x)$, $x\in X$ is uniformly bounded in $x$ by a finite number), and especially for finitely generated groups equipped with the Cayley distance. As target spaces $Y$, one is typically interested in ``nice'' Banach spaces (the Hilbert space $\ell^2$, the spaces $\ell^p$, $p>1$, etc.) 

In view of the example of Dranishnikov {\em et al.}, the standard ultraproduct technique implies immediately that there exists a locally finite metric space (that is, every ball of finite radius contains finitely many points) non-embeddable into $\ell^2$ \cite{dgly}. The same question for spaces of bounded geometry is more difficult. It was resolved by Gromov \cite{gromov00}, who has noticed that a disjoint sum of graphs forming an expander family and equipped with the path distance gives such an example of a metric space of bounded geometry. (This construction is treated in detail e.g. in the book \cite{roe}.) In the same paper \cite{gromov00}, Gromov has outlined a probabilistic construction of a finitely generated infinite group into which a disjoint sum of graphs as above can be coarsely embedded and which therefore admits no coarse embedding into $\ell^2$. 

Recent results by Kasparov and Yu \cite{ky} have brought interest to a more general version of the same coarse embedding problem, where $\ell^2$ is replaced with a uniformly convex Banach space. For some remarks on the problem, see \cite{lw}.

Recall that a Banach space $E$ is {\em uniformly convex} if the function $\delta\colon (0,2]\to \R$ (the {\em modulus of convexity} for $E$), defined by
\[\delta(\e)=\inf\left\{1-\norm{\frac{x+y}2}\colon x,y\in E,~~ \norm x,\norm y\leq 1,~~\norm{x-y}\geq\e\right\},\]
is strictly positive for all $\e\in (0,2]$. 

The spaces $\ell^p$, $1<p<\infty$, are uniformly convex. The Hilbert space $\ell^2$ coarsely embeds into $\ell_p$ for $1\leq p\le\infty$, and the spaces $\ell^p$, $1<p\leq 2$, admit coarse embeddings into each other (Nowak \cite{nowak05,nowak06}). At the same time, $\ell^p$ do not coarsely embed into $\ell^2$ for $p>2$ (W.B. Johnson and  Randrianarivony \cite{jr}). (Notice that all the proofs mentioned in this paragraph are also emulations of analogous results obtained for uniform embeddings much earlier.)

While working on this article, the author became aware of an announcement by Vincent Lafforgue \cite{lafforgue} that a family of expander graphs associated to the groups ${\mathrm{SL}}(3)$ over $p$-adic fields does not embed into any uniformly convex Banach space. 

Another problem that has stimulated the present investigation is a question of the existence of a metric space that admits no uniform embedding into a {\em reflexive} Banach space. (Cf. e.g. question 6.5 in \cite{megrelishvili07}.) This question was answered in the negative by Kalton \cite{kalton}.

In the present paper, we show that the universal Urysohn metric space $\Ur$ admits neither uniform nor coarse embedding into a uniformly convex Banach space. Of course, the space $\Ur$ is far from having a bounded geometry, quite on the contrary. Besides, the results turn out to be outdone not only by Lafforgue's result mentioned above, but by a remarkable theorem of Kalton \cite{kalton}: the space $c_0$ admits neither uniform nor coarse embedding into a {\em reflexive} Banach space. (I have learned about the yet unpublished paper of Kalton from the arXiv preprint \cite{BaL}, which was submitted after the first version of the present article.) 

However, our method of proof is rather different from the existing methods for showing non coarse embeddability of various metric spaces. We are using the Hrushovski--Solecki--Vershik property of the Urysohn space as a basis for an averaging argument, which could be of interest on its own. This is why the present author feels the publication of this work still has some merit. 

\section{The Urysohn metric space}

The universal Urysohn metric space $\Ur$ is determined uniquely, up to an isometry, by the following description: $\Ur$ is a complete separable metric space which contains an isometric copy of every separable metric space and is {\em ultrahomogeneous,} that is, every isometry between two finite subspaces of $\Ur$ can be extended to a global self-isometry of the space $\Ur$. 

Chapter 5 in the author's book \cite{P06} is an introduction to the Urysohn space rather well suited for our needs. Another highly-recommended introductory source is Melleray's article \cite{melleray} written for the present volume. 
Other self-contained introductions to the Urysohn space can be found in Gromov's book \cite{Gr} and articles by Vershik \cite{Ver98,Ver04} and Uspenskij  \cite{Usp98}.

One obtains numerous variations of the Urysohn space by restricting the set of all possible values that the distances in the metric spaces appearing in the above description of $\Ur$ can assume. Such variations include the {\em integer-valued} Urysohn space $\Ur_{\Z}$ and the {\em rational Urysohn space,} $\Ur_{\Q}$. One can also consider the {\em bounded Urysohn space of diameter one}, $\Ur_1$, where all distances between two different points take values in the interval $[0,1]$. It is an easy exercise, to show that $\Ur_1$ is isometric to the sphere of radius $1/2$ taken in $\Ur$ around any point. The space $\Ur_1$ is also known as the {\em Urysohn sphere.}

Recall that Rado's {\em infinite random graph,} $R$, is defined by the following conditions: it is a simple unoriented graph that is universal for all countable graphs (that is, contains an isomorphic copy of every countable graph as an induced subgraph)
and ultrahomogeneous in the sense that every isomorphism between two finite induced subgraphs extends to a global automorphism of $R$ \cite{cameron}. The Urysohn space $\Ur_{0,1,2}$ whose distances take only values $0$, $1$ and $2$, is easily shown to be isometric to the infinite random graph $R$ equipped with the path-distance (where every edge has length $1$). 

Now we need to recall another notion, due to Uspenskij \cite{Usp98}. Say that a subspace $X$ of a metric space $Y$ is {\em $g$-embedded} into $Y$ if there is a  group homomorphism $h\colon \Iso(X)\to\Iso(Y)$, continuous with regard to the topology of simple convergence on both groups and such that for every isometry $i\in\Iso(X)$ one has
 \[h(i)\vert_X=i.\]
In other words, one can simultaneously extend all isometries from $X$ to $Y$, preserving the algebraic operations and in a ``continuous way''.
 
One of the most useful results of the theory of the Urysohn space, due to Uspenskij \cite{Usp90}, says that every separable metric space can be $g$-embedded into the Urysohn space $\Ur$. See also \cite{Usp98,P06,melleray}. The same method (Kat\v etov extensions) lead to the following observation, also due to Uspenskij.

\begin{theorem} Every compact metric subspace $X$ of $\Ur$ is $g$-embedded. \qed
\end{theorem}

The following result is obtained along the same lines.

\begin{theorem}
Every finite metric subspace of the integer Urysohn space $\Ur_\Z$ is $g$-embedded.
\end{theorem}

\begin{pf} The proof is modelled on Uspenskij's original argument, as presented, e.g., in our book \cite{P06} on pp. 109-111, 114-115. A function $f\colon X\to\R$ on a metric space is called a {\em Kat\v etov function} if 
\begin{equation}
\vert f(x)-f(y)\vert\leq d_X(x,y)\leq f(x) + f(y)
\label{eq:flood}
\end{equation}
for all $x,y\in X$. Such functions are exactly distances from points in metric extensions of $X$. A 1-Lipschitz function $f$ as above is
{\em controlled by a metric subspace} $A\subseteq X$
if $f$ is the largest 1-Lipschitz function on $X$
having a given restriction to $A$. The set of all Kat\v etov functions controlled by finite subspaces of $X$, equipped with the supremum metric, is denoted $E(X)$, and contains $X$ in a canonical way under the Kuratowski embedding associating to an $x\in X$ the distance function from $x$. The completion of the union of a recursively built chain of $n$-fold extensions $E^n(X)=E(E^{n-1}(X))$ is isometrically isomorphic to the Urysohn space $\Ur$.

Denote by $E_\Z(X)$ the collection of all integer-valued Kat\v etov functions on $X$ controlled by finite subsets, equipped with the uniform distance $d(f,g)=\sup\{\abs{f(x)-g(x)}\colon x\in X\}$. This metric space is nontrivial, because it contains, for instance, all distance functions $d(x,-)$ from points in $X$. Moreover, the Kuratowski embedding
\[X\ni x\mapsto d(x,-)\in E_\Z(X)\]
is easily seen to be an isometric embedding. We will thus identify $X$ with a metric subspace of $E_\Z(X)$. The union of an infinite chain of iterated $n$-fold Kat\v etov extensions of $X$ can be verified to be isometric to $\Ur_\Z$ (cf. \cite{lanvt}). If $G$ is a group acting on $X$ by isometries, then the action lifts to the Kat\v etov extension $E_\Z(X)$ through the left regular representation, and the Kuratowski embedding is $G$-equivariant under this lifting. This means that the action of $G$ goes all the way through to $\Ur$. \qed
\end{pf}

The space $E_\Z(X)$ also appears in \cite{lanvt} and \cite{Ver05b}. Here is another application of this construction.

From the viewpoint of coarse geometry, all bounded metric spaces, such as $\Ur_1$ or $R$, are alike: they are coarsely equivalent to a single point. 
Two metric spaces $X$ and $Y$ are {\em coarsely equivalent} if there is a pair of maps $f\colon X\to Y$, $g\colon Y\to X$ such that the compositions $fg$, $gf$ are uniformly close to the corresponding identity maps, that is, the functions $d_X(gf(x),x)$ and $d_Y(fg(y),y)$ are bounded on $X$ and $Y$, respectively.

The real line $\R$ is coarsely equivalent to the subspace $\Z$. As a function $f\colon\Z\to\R$, one can take the canonical embedding, and as a function $g\colon\R\to\Z$, the integer part function. 
This observation generalizes as follows: every metric space is coarsely equivalent to its subspace forming an $\e$-net for some $\e>0$. As one of the authors of \cite{lanvt} (Lionel Nguyen Van Th\'e) has pointed out to the present author, the proof of Proposition 1 (Section 2.1 in \cite{lanvt}) can be modified so as to establish the following result.

\begin{theorem}[Jordi L\'opez-Abad and Lionel Nguyen Van Th\'e] The Urysohn space $\Ur$ contains as a 1-net an isometric copy of the integer Urysohn space $\Ur_{\Z}$. \qed
\label{th:jordilionel}
\end{theorem}

\begin{corollary}
\label{c:jordilionel}
The Urysohn space $\Ur$ is coarsely equivalent to the integer-valued Urysohn space $\Ur_\Z$. \qed
\end{corollary}

Since the composition of two coarse embeddings is a coarse embedding, and every coarse equivalence is a coarse embedding, it follows from Corollary \ref{c:jordilionel} that for the purpose of considering coarse embeddings, there is no difference between $\Ur$ and $\Ur_{\Z}$. 

\section{\label{s:hsv}
The Hrushovski--Solecki--Vershik homogeneity property}

\begin{definition}
Let us say that a metric space $X$ has the {\em Hrushovski--Solecki--Vershik property} if for every finite subspace $Y$ of $X$ there exists a finite $Z$ with $Y\subseteq Z\subseteq X$ such that every partial isometry of $Y$ extends to a self-isometry of $Z$.
\label{d:hsv}
\end{definition}

\begin{proposition}
\label{p:implies}
For a complete separable metric space $X$ the Hrushovski--Solecki--Vershik property implies ultrahomogeneity.
\end{proposition}

\begin{pf}
Required global isometries of $X$ are built up recursively. Namely, if $A$ and $B$ are two finite metric subspaces of $X$, and $i\colon A\to B$ is an isometry, then there exists a finite $Y\subseteq X$ and a self-isometry $j_0$ of $Y$ whose restriction to $A$ coincides with $i$. Now enumerate an everywhere dense subset $X^\prime=\{x_i\colon i\in\N_+\}$ of $X$, and choose an increasing chain 
\[Y=Y_0\subseteq Y_1\subseteq\ldots\subseteq Y_n\subseteq\ldots\]
of subspaces of $X$ and their self-isometries $j_n$ in such a way that $Y_{n+1}\supseteq Y_n\cup\{x_n\}$ and $j_{n+1}\vert_{Y_n}=j_n$. The mapping $j$ defined by the rule $j(x)=j_n(x)$ whenever $x\in Y_n$ is a self-isometry of $X^\prime$, and it extends by continuity over all of $X$. 
\end{pf}

\begin{corollary}
A metric space $X$ has the Hrushovski--Solecki--Vershik property if and only if for every finite subspace $Y$ of $X$ there exists a finite $Z$ with $Y\subseteq Z\subseteq X$ such that every partial isometry of $Y$ extends to a global self-isometry of $X$ taking $Z$ to itself. \qed
\end{corollary}

The converse to Proposition \ref{p:implies} is not true: the Hrushovski-Solecki-Vershik property does not in general follow from ultrahomogeneity.

The topology of pointwise convergence on the isometry group $\Iso(X)$ of a metric space $X$ coincides with the compact-open topology and is a group topology. If $X$ is in addition separable and complete, then $\Iso(X)$ is a Polish group. The following belongs to Vershik.

\begin{proposition}
If a separable metric space $X$ has the Hrushovski-Solecki-Vershik property, then the isometry group $\Iso(X)$ contains an everywhere dense locally finite subgroup.
\label{p:dlf}
\end{proposition}

\begin{pf} A dense locally finite subgroup is built recursively, much like in the proof of Proposition \ref{p:implies}.
\end{pf}

\begin{example}
The Hilbert space $\ell^2$, as well as all Euclidean spaces $\ell^2(n)$ (including the real line $\R$), are ultrahomogeneous (\cite{Blum}, Ch. IV, \S 38, Property 4). However, each of them
fails the Hrushovski--Solecki--Vershik property. 
\end{example}

Indeed,
let $x\neq 0$ be an element of a Euclidean space $E$ (finite or infinite-dimensional), and let $T_x$ denote the translation $y\mapsto y+x$. Let $Y=\{0,x,2x,3x\}$, $A=\{0,x,2x\}$, and let $i=T_x\vert_{A}$.
Let $j$ be an isometry of some larger subspace $Z\supseteq Y$ of $E$ extending the partial isometry $i$. The points $j(0)=x$, $j(x)=2x$ and $j(2x)$ form a metric space isometric to $(0,x,2x)$ and therefore lie on a straight line, and consequently one must have $j(2x)=3x$. An inductive argument shows that $Z$ contains all elements of the form $nx$ and thus is infinite. \tri

A similar argument (using a rotation by an irrational angle along a grand circle instead of a translation along a straight line) gives the following.

\begin{example}
The unit sphere $\s$ of every Hilbert space of positive dimension (including the circle $\s^1$ and the infinite-dimensional unit sphere $\s^\infty$ of $\ell^2$) does not have the Hrushovski--Solecki--Vershik property.
\end{example}

The sphere $\s^\infty$ has a somewhat weaker property: every finite collection of isometries can be simultaneously approximated in the strong operator topology with elements of a finite subgroup of isometries (Kechris \cite{kechris}, a remark on page 186). Notice that the sphere $\s^\infty$ is ultrahomogeneous.

The following two results explain the origin of the name that we gave to the property above. The first one is an equivalent reformulation of a result by Hrushovski \cite{hrushovski}.

\begin{theorem} 
\label{th:h}
The infinite random graph (or, equivalently, the Urysohn space $\Ur_{0,1,2}$) has the Hrushovski--Solecki--Vershik property. \qed
\end{theorem} 

The second result has been established by Solecki \cite{solecki} and, independently, by Vershik \cite{Ver05b}. The following is an equivalent reformulation.

\begin{theorem}
\label{th:sv}
The Urysohn metric space $\Ur$, the rational Urysohn space $\Ur_\Q$ and the integer Urysohn space $\Ur_{\Z}$ all have the Hrushovski--Solecki--Vershik property. \qed
\end{theorem} 

Solecki's proof is a corollary of deep results of the paper \cite{HL} of Herwig and Lascar, while Vershik's proof is direct. In the Appendix we present a deduction of the theorem from results of \cite{HL}, although in a somewhat different way from Solecki's argument.

Hrushovski's theorem admits a very simple direct combinatorial proof, cf. \cite{HL}, Sec. 4.1. This author has been unable to find an analogue for finite metric spaces.
% 
% \begin{remark} 
% \label{r:stronger}
% Let us say provisionally that a metric space $X$ satisfies the {\em diameter-preserving Hrushovski--Solecki--Vershik property} if the finite metric space $Z$ selected as in Definition \ref{d:hsv} can be selected so as to have $\diam Z=\diam Y$.
% 
% The above proof establishes that the metric space $\Ur_S$ has the diameter-preserving Hrushovski--Solecki--Vershik property.
% \end{remark}
% 
% Indeed, one can replace the set $S$, without any loss in generality, with the convex subset $S\cap [0,\diam X]$ of $T$. \qed

\section{Averaging distances}

Let $X_\alpha,\alpha\in A$ be a family of metric spaces. Choose a point $x_\alpha$ in each $X_\alpha$.
Define the set
\[\ell^\infty(X_\alpha,x_{\alpha},A)= \left\{y\in \prod_{\alpha\in A}X_\alpha\mid \sup_{\alpha\in A}d(x_\alpha,y_\alpha)<\infty\right\}.\]
Now let $\xi$ be an ultrafilter on the index set $A$. Equip $\ell^\infty(X_\alpha,x_{\alpha},A)$ with the pseudometric
\[d(y,z)=\lim_{\alpha\to\xi}d(y_\alpha,z_{\alpha}).\]
The {\em metric space ultraproduct} along $\xi$ of the family $(X_\alpha)$ centred at $(x_\alpha)$ is the metric quotient of the pseudometric space $\ell^\infty(X_\alpha,x_{\alpha},A)$. If $\xi$ is a free ultrafilter, then the ultraproduct is complete.

In a particular case where every $X_\alpha$ is a normed space and $x_\alpha=0$, we obtain the familiar concept of the ultraproduct of a family of normed spaces.
% 
% Let $E_\alpha,\alpha\in A$ be such a family, and let
% $\xi$ be an ultrafilter on the index set $A$. The {\em Banach space ultraproduct} of
% the family $(E_\alpha)$ along the ultrafilter $\xi$ is the linear space
% quotient of the $\ell^\infty$-type direct sum
% $E=\oplus^{\ell^\infty}_{\alpha\in A}E_\alpha$ by the ideal ${\mathcal I}_\xi$
% formed by all collections $(x_\alpha)_{\alpha\in A}\in E$ having the property
% \[\lim_{\alpha\to\xi}x_{\alpha} =0.\]
% The ultraproduct is equipped with the norm
% \[\norm x = \lim_{\alpha\to\xi}x_{\alpha},\]
% where $(x_\alpha)$ is a representative of the equivalence class $x$.
% If the ultrafilter $\xi$ is free, the ultraproduct is a Banach space.
In nonstandard analysis, Banach space ultraproducts are known as \textit{nonstandard hulls}. For more on them, see \cite{HI} and references therein.

However, even in the case where $x_\alpha\neq 0$, the construction is of interest.

\begin{proposition}
The ultraproduct of a family of normed spaces $X_\alpha$ centred at an arbitrary family of points $(x_\alpha)$ becomes an affine normed space in a natural way, and for different choices of $(x_\alpha)$ all affine normed spaces arising in this way are pairwise affinely isomorphic and isometric.
\end{proposition}

\begin{pf} 
Every linear translation in the product $\prod_{\alpha\in A}X_\alpha$ preserves the values of the distance $\lim_{\alpha\to\xi}d(y_\alpha,z_{\alpha})$ and consequently defines an isometry between ultraproducts centered at different points. Since every self-isometry of a normed space is an affine map, it follows that every isometry between two such ultraproducts transports the affine structure in a unique way. An altervative way to furnish each ultraproduct with an affine structure is through noticing that the space $\ell^\infty(X_\alpha,x_{\alpha},A)$ is an affine subspace of the linear space $\prod_{\alpha\in A}X_\alpha$, of which the ultraproduct is a quotient affine space.
\end{pf}

For a normed space $E$ and a set $Z$, denote by $\ell^2(E,Z)$ the $\ell^2$-type sum of $\abs Z$ copies of $E$.

\begin{theorem}
Let $G$ be a locally finite group acting by isometries on a metric space $X$ and having a dense orbit. Suppose that $X$ admits a mapping $\phi$ into a normed space $E$ such that for some functions $\rho_1,\rho_2\colon\R_+\to\R_+$ one has
\begin{equation}
\rho_1(d_X(x,y))\leq \norm{\phi(x)-\phi(y)}\leq \rho_2(d_X(x,y)).
\label{eq:embe}
\end{equation}
Then there is a map $\psi$ of $X$ into a Banach space ultrapower of  $\ell^2(E,\xi)$, satisfying the same inequalities (\ref{eq:embe}) and such that
the action of $G$ on $\psi(X)$ extends to an action of $G$ by affine isometries on the affine span of $\psi(X)$.
% If in addition the action of $G$ on $X$ has bounded orbits, then one can assume the above affine action to be a linear representation by isometries.
% the image $\psi(X)$ is a metric transform of $X$: the distance $\norm{\psi(x)-\psi(y)}$ is a function of $d(x,y)$ alone. In addition, 
\end{theorem}

\begin{pf}
Denote by $\Xi$ the set $\Xi$ of all finite subgroups of $G$.
For every $F\in\Xi$ define a map $\psi_F\colon X\to \ell^2(E,F)$ as follows:
\[\psi_F(x)(g)=\frac{1}{\abs{F}^{1/2}}\psi(g^{-1}x).\]
Since the distance $d=d_X$ is $G$-invariant, the map $\psi_F$ satisfies the inequalities (\ref{eq:embe}), e.g. for the right hand side inequality:
\begin{eqnarray*}
\norm{\psi_F(x)-\psi_F(y)}&=& 
\left( \sum_{g\in F}\frac{1}{\abs{F}}\norm{\psi(g^{-1}x)-\psi(g^{-1}y)}_E^2 \right)^{1/2}\\
&\leq& \frac{1}{\abs{F}^{1/2}} \left( \abs{F} \cdot\rho_2^2(d(x,y))\right)^{1/2}\\
&=& \rho_2(d(x,y)).
\end{eqnarray*}
Make $\ell^2(E,F)$ into an $F$-module via the left regular representation:
\[^hf(g)=f(h^{-1}g),~~g,h\in F.\]
The map $\psi_F$ is $F$-equivariant: for every $g\in F$,
\begin{eqnarray*}
^h\psi_F(x)(g)&=& \psi_F(x)(h^{-1}g)\\
&=& \frac{1}{\abs{F}^{1/2}}\psi(g^{-1}hx)\\
&=&\psi_F(hx)(g).
\end{eqnarray*}

Choose an ultrafilter $\xi$ on $\Xi$ with the property that for each $F\in\Xi$ the set of all $\Phi\in\Xi$ containing $F$ is in $\xi$. Select a point $x^\ast\in X$ whose $G$-orbit is dense in $X$. Denote by $V$ the ultraproduct of the spaces $\ell^2(E,F)$, $F\in\Xi$, along the ultrafilter $\xi$, centred at the family of points $(\psi_F(x^\ast))_{F\in\Xi}$. 

For every $h\in G$,
the family $(\psi_F(hx^\ast))_{F\in\Xi}$ is at a finite distance from $(\psi_F(x^\ast))_{F\in\Xi}$: 
\begin{eqnarray*}
\sup_{F\in\Xi}
\norm{\psi_F(hx^\ast))-\psi_F(x^\ast))} &=& 
\sup_{F\in\Xi}
\frac{1}{\abs{F}^{1/2}}\left(\sum_{g\in F}\norm{\psi(g^{-1}hx^\ast)-\psi(g^{-1}x^\ast)}_E^2
\right)^{1/2} \\
&\leq& \sup_{F\in\Xi}\frac{1}{\abs{F}^{1/2}}\left(\abs{F}\cdot \rho_2(d(g^{-1}hx^\ast,g^{-1}x^\ast)^2\right)^{1/2}\\
&=&\rho_2(d(hx^\ast,x^\ast)).
\end{eqnarray*}
This has two consequences. First, since the $G$-orbit of $x^\ast$ is dense in $X$, the family of mappings $(\psi_F)$ determines a well-defined mapping $\psi\colon X\to V$ (that is, for every $x\in X$ the image $\psi(x)$ is a well-defined element of $V$). Second, for every $h$ the translation by $h$ determines an isometry of $V$, and in this way the group $G$ acts on $V$ by isomeries.
making the mapping $\psi$ $G$-equivariant. 

It remains to notice that the space $\ell^2(E,\xi)$ contains all spaces $\ell^2(E,F)$ as normed subspaces, so a suitable metric space ultrapower of $\ell^2(E,\xi)$ contains our metric space ultraproduct of $\ell^2(E,F)$. Finally, each metric ultrapower of $\ell^2(E,\xi)$ is isometrically affinely isomorphic to the Banach space ultrapower.
% 
% If in addition the action of $G$ on $X$ has bounded orbits, then the family $(\psi_F(x^\ast))_{F\in\Xi}$ is at a finite distance from zero element, and it is clear that the translations by $g$ (which are all linear in the spaces $\ell^2(E,F)$) determine a linear isometric transformations of the corresponding Banach space ultrapower. 
\end{pf}

Let us introduce a natural concept: say that a group $G$ of isometries of a metric space $X$ is {\em almost $n$-transitive} if for every $\e>0$ and every isometry $i$ between two subspaces $A,B\subseteq X$ of cardinality $n$ each there is a $g\in G$ with the property that for all $a\in A$ one has $d(i(a),ga)<\e$. 

\begin{remark}
Every dense locally finite subgroup of the group of isometries $\Iso(X)$ of a separable metric space having the Hrushovski--Solecki--Vershik property (Proposition \ref{p:dlf}) is almost $n$-transitive for every $n$. 
\label{r:dense}
\end{remark}

\begin{corollary} 
Let $X$ be a metric space admitting a locally finite almost $2$-transitive group of isometries. 
% the Hrushovski--Solecki--Vershik property, and 
Let $\psi$ be a mapping of $X$ into a normed space $E$ satisfying inequalities (\ref{eq:embe}).
% 
% such that for some functions $\rho_1,\rho_2\colon\R_+\to\R_+$
% \begin{equation}
% \rho_1(d_X(x,y))\leq \norm{\phi(x)-\phi(y)}\leq \rho_2(d_X(x,y)).
% \label{eq:embe1}
% \end{equation}
Then there is a map $\psi$ of $X$ into a Banach space ultrapower of  $\ell^2(E,\xi)$, satisfying the same inequalities (\ref{eq:embe}) and such that the image $\psi(X)$ is a metric transform of $X$: the distance $\norm{\psi(x)-\psi(y)}$ is a function of $d(x,y)$ alone. Furthermore, the action of the group of isometries of $X$ extends to a representation by affine isometries on the affine span of $\psi(X)$, making $\psi$ an equivariant map.
\label{l:main}
\qed
\end{corollary}

\section{Non-existence of uniform embeddings}

Uniform convexity is a {\em metric} property of a Banach space, which can be lost if a norm is replaced by an equivalent one. The corresponding property of Banach spaces invariant under isomorphisms is {\em superreflexivity.} A Banach space is superreflexive if every Banach space ultrapower of it is a reflexive Banach space. It can be shown that a Banach space is superreflexive if and only if it admits an equivalent uniformly convex norm (cf. \cite{fhhmpz}). In our context, speaking of superreflexive Banach spaces is more appropriate, because both coarse and uniform structures are invariant under Banach space isomorphisms.

\begin{theorem} The universal Urysohn metric space $\Ur$ cannot be embedded, as a uniform subspace, into a superreflexive Banach space.
\label{th:non-uniform}
\end{theorem}

\begin{pf}
Assume, towards a contradiction, that an embedding, $\phi\colon \Ur\to E$, does exist. Choose a dense locally finite subgroup of $\Iso(\Ur)$ (Proposition \ref{p:dlf}). This group is almost $n$-transitive for each $n$. By Corollary \ref{l:main}, there exist a mapping, $\psi$, of $\Ur$  into the Banach space ultrapower $V$ of $\ell^2(E)$, as well as a representation $\pi$ of $G$ by affine isometries of the affine span $F$ of $\psi(\Ur)$ in $V$, making $\psi$ equivariant. 

Since the mapping $\psi$ is a uniform isomorphism on its image and in particular a homeomorphism, the topology of pointwise convergence on $X$ on the group $G$ coincides with the topology of pointwise convergence on $\psi(X)$ and consequently on the affine span of $\psi(X)$. The action of $G$ on $F$ extends by continuity to an action of all of $\Iso(\Ur)$ on $F$, and $\psi$ still remains equivariant. Consequently, the affine representation of $\Iso(\Ur)$ on $F$ is faithful. 

This representation is a continuous homomorphism from $G$ to the group of affine isometries of $F$, that is, the semidirect product $\Iso(F)\ltimes F_+$. Here $\Iso(F)$ denotes the group of linear isometries, while $F_+$ is the additive group of $F$ upon which $\Iso(F)$ acts in a canonical way. Let $\pi\colon \Iso(F)\ltimes F_+\to \Iso(F)$ be the standard projection (a quotient homomorphism). 

The Polish group $\Iso(\Ur)$ is universal \cite{Usp98}. In particular, it contains, as a topological subgroup, the group $\Homeo_+[0,1]$ of all homeomorphisms of the unit interval preserving endpoints, equipped with the standard compact-open topology. We have, therefore, a faitful continuous affine representation of this group in a superreflexive Banach space.

According to Megrelishvili \cite{megrelishvili01}, the only continuous representation of $\Homeo_+[0,1]$ by linear isometries in a reflexive Banach space is the trivial (identity) representation. Therefore, the composition of three homomorphisms
\[\Homeo_+[0,1]\hookrightarrow \Iso(\Ur)\to \Iso(F)\ltimes F_+\overset{\pi}{\to} \Iso(F)\]
is a trivial map. This means the image of $\Homeo_+[0,1]$ is contained, as a topological subgroup, in the kernel of $\pi$, that is, the abelian Polish group $F_+$, which is absurd.
\qed
% 
% Taking into account that the group $\Iso(\Ur)$ is topologically simple \cite{Usp98}, that is, contains no proper closed normal subgroups, it is easy to conclude that $\Iso(\Ur_1)$ admits no non-trivial strongly continuous representations by isometries in reflexive Banach spaces. To obtain a contradiction, notice that $F$ is reflexive, and contains at least one non-trivial orbit of the representation $\pi\vert_{\Iso(\Ur_1)}$, namely a uniformly isomorphic copy of the Urysohn sphere $\Ur_1$. \qed
\end{pf}

\section{Non-existence of coarse embeddings}

We will begin by recalling a useful test for a space not to be (super)reflexive. Let $\e>0$ and let $n$ be either a natural number or the symbol $\infty$. An $(n,\e)${\em -tree} in a normed space $E$ is a binary tree $T$ of depth $n$ whose nodes are elements of $E$ such that for every node $x$ its children nodes $y$ and $z$ have the properties: $x=(y+z)/2$ and $\norm{y-z}\geq\e$. 

\begin{theorem}[See e.g. \cite{fhhmpz}, p. 295] If a normed space $E$ contains a uniformly bounded family of $(n,\e)$-trees for some $\e>0$ and all natural $n$, then $E$ is not superreflexive. If $E$ contains a bounded $(\infty,\e)$-tree for some $\e>0$, then $E$ is non-reflexive. \qed
\end{theorem}

Here is a consequence that we will be using. 

\begin{corollary}[Cf. a similar statement in \cite{fhhmpz}, Exercise 9.22, p. 308]
Let a normed space $E$ have the following property: for some $\e>0$ and $M>0$ and for every $n$, the $M$-ball around zero contains a sequence of closed convex subsets $K_i$, $i=1,2,\ldots,n$ such that $K_{2i}\cup K_{2i+1}\subseteq K_i$ and 
$K_{2i}$ and $K_{2i+1}$ are at a distance at least $\e$ from each other for all $i$. Then $E$ is not superreflexive.
\label{c:corol}
\end{corollary}

\begin{pf} For every $n$ one can easily construct a $(n,\e)$-tree contained in the $M$-ball of $E$ by recursion, starting with the leaves and using the fact that the midpoint of two nodes belonging to $K_{2i}$ and $K_{2i+1}$, respectively, is contained in $K_i$. The assumption on the distance between $K_{2i}$ and $K_{2i+1}$ assures that the two children nodes are always at least $\e$-apart from each other. \qed
\end{pf} 

We need to recall a classical result by Day.

\begin{theorem}[Day \cite{day}]
The $\ell_p$-type direct sum of normed spaces, $1<p<\infty$, is uniformly convex if and only if they have a common modulus of convexity. \qed
\end{theorem}

Recall that a normed space $E$ is {\em uniformly smooth} (\cite{BL}, Appendix A; \cite{fhhmpz}, Ch. 9) if $\norm{x+y}+\norm{x-y}=o(\norm y)$ as $\norm y\to 0$, uniformly for all $x$ in the unit sphere of $E$. More precisely, the {\em modulus of uniform smoothness}, $\delta$, of $E$ is defined for $\tau>0$ by
\[\rho(\tau)=\sup\left\{\frac{\norm{x+\tau h}+\norm{x-\tau h}}2\colon \norm x =\norm h =1\right\}.\]
Now $E$ is uniformly smooth if and only if 
\[\lim_{\tau\to 0}\frac{\rho(\tau)}{\tau}=0.\]

An ultraproduct of a family of uniformly smooth normed spaces of the same modulus of smoothness is again uniformly smooth, which is straightforward. The same is true of the $\ell^p$-direct sum, $1<p<\infty$. (The result can be deduced from Day's theorem by using duality, as there is a correspondence between the moduli of uniform convexity and of uniform smoothness of $E$ and of $E^\ast$.) 

Every superreflexive space admits an equivalent norm that is both uniformly convex and uniformly smooth. A space $E$ is uniformly convex iff $E^\ast$ is uniformly smooth, and vice versa. (See \cite{fhhmpz}, Ch. 9.)

Every point $x$ of the unit sphere of a uniformly smooth normed space $E$ is a smooth point, that is, there exists a unique linear functional $j(x)=\varphi\in E^\ast$ of norm one such that $\varphi(x)=1$. (See e.g. \cite{BL}, p. 70.)

\begin{theorem}
The universal Urysohn metric space $\Ur$, as well as the integer-valued Urysohn space, $\Ur_\Z$, do not admit a coarse embedding into a superreflexive Banach space.
\end{theorem}

\begin{pf}
Assume, towards a contradiction, that a coarse embedding $\phi\colon\Ur_\Z\to E$ exists. 
Since the space $E$ is superreflexive, it admits an equivalent norm that is at the same time uniformly convex and uniformly smooth. 

% Choose an arbitrary point $x_0\in\Ur_\Z$. 
By Corollary \ref{l:main}, there exist a coarse embedding $\psi$ of $\Ur_\Z$ into an ultrapower of $\ell^2(E)$ and a strongly continuous representation of $\Iso(\Ur_\Z)$ by affine isometries in a closed subspace $F$ spanned by $\psi(\Ur_\Z)$, such that $\psi$ is an $\Iso(\Ur_\Z)$-equivariant mapping. 
Since $\ell^2$-sums and ultraproducts preserve uniform convexity and uniform smoothness, and these properties are inherited by normed subspaces, the norm on the space $F$ is both uniformly convex and uniformly smooth. 

It follows from the same Corollary \ref{l:main} that the coarse embedding $\psi$
is a metric transform:
% that is, $\norm{\psi(x)-\psi(y)}$ only depends on $d(x,y)$ for all $x,y\in \Ur_{\Z}$:
\[\forall x,y\in\Ur_{\Z},~~\norm{\psi(x)-\psi(y)}=\rho(d(x,y)),\]
where $\rho(r)\to\infty$ as $r\to\infty$. By rescaling the norm of $F$ if necessary, one can assume without loss in generality that $\rho(1)=1$.

Choose an arbitrary point $z_0\in\Ur_\Z$ and denote $G=\Iso(\Ur_{\Z})_{z_0}$ the isotropy subgroup of $z_0$ in $\Iso(\Ur_{\Z})$. Make the affine space $F$ into a linear space by setting $0_F=\psi(z_0)$. The action of $G$ on $F$ becomes a linear continuous isometric representation.

Fix a natural number $m$ with $\rho(m)\geq 6$, and 
% One can assume without loss in generality that $\rho_1(1)=1$. Indeed, for some $n\in\N_+$ one has $\rho_1(n)>0$. The version of the Urysohn space $\Ur_{n\N}$, whose distances take values in $n\N$, isomerically embeds into $\Ur_{\Z}$, and by scaling all distances in $\Ur_{n\N}$ by a factor of $1/n$, we get an isomeric copy of $\Ur_\Z$. Now it remains to rescale the norm in the target Banach space by a factor of $\rho_1(n)^{-1}$. 
% 
% Since the function $\rho_1$ is unbounded, there is an even natural number $m$ with $\rho_1(m)>6\rho_2(1)$. The restriction of $\psi$ to the sphere $S_m(x_0)$ is a metric transform, that is, there exists a function $\varsigma\colon\R\to\R$ such that for all $x,y\in S_m(x_0)$ one has
% \[\norm{\psi(x)-\psi(y)}=\varsigma(d(x,y)).\]
% In particular, if $d(x,y)=1$, then 
% \[1\leq\norm{\psi(x)-\psi(y)}=\varsigma(1)\leq\rho_2(1),\]
% and so
% \begin{equation}
% \varsigma(m)>6\varsigma(1).
% \label{eq:6}
% \end{equation}
let $x,y\in S_m(z_0)$ be such that $d(x,y)=m$. Then $\norm{\psi(x)-\psi(y)}=\rho(m)$ and by the Hahn--Banach theorem, there is a linear functional $\varphi\in F^\ast$ of norm $1$ such that 
\[\varphi(\psi(x)-\psi(y))=\rho(m).\]
Since $F$ is uniformly smooth, such a $\varphi$ is unique (the support functional of $\psi(x)-\psi(y)$).

Let $x=x_0,x_1,\ldots,x_{m-1},x_m=y$ be a sequence of elements in $S_m(z_0)$ satisfying $d(x_i,x_{j})=\abs{i-j}$ for all $i,j=0,1,\ldots,m$. Any two subsequent values $\varphi(\psi(x_i))$, $i=0,1,\ldots,m$,
differ by at most $1$. 

There is an isometry $f\in G$ interchanging $x_j$ and $x_{n-j}$ for every $j$. In particular, $f$ flips $x$ and $y$. The corresponding linear isometry of $E^\ast$ will take $\varphi$ to the support functional of $\psi(y)-\psi(x)$, that is, to $-\varphi$. This means that $\varphi(\psi(x_j))=-\varphi(\psi(x_{m-j}))$ for every $j=0,1,2,\ldots,m$, in particular, $\varphi(\psi(x))=-\varphi(\psi(y))$. Without loss in generality one can assume that $\varphi(\psi(x))$ is negative. 

Denote $k=\max\{j=0,1,2\ldots,m\colon \varphi(\psi(x_j))<0\}$. Then $\varphi(\psi(x_{k+1}))\geq 0$. 
Let $z$ and $w$ be arbitrary points of $S_{m}(z_0)$ such that $d(x,z)=k$, $d(x,w)=k+1$, $d(y,z)=m-k$, and $d(y,w)=m-k-1$. There exists a $g\in G$ stabilizing $x$ and $y$ and taking $z\mapsto x_k$ and $w\mapsto x_{k+1}$. The extension of $g$ to $F^\ast$ will leave $\varphi$ fixed, because of its uniqueness as a support functional. We conclude that $\varphi(\psi(z))=\varphi(\psi(x_k))$ and $\varphi(\psi(w))=\varphi(\psi(x_{k+1}))$. This can be summarized as follows: the functional $\varphi$ assumes the constant value $\varphi(\psi(x_k))<0$ at all points of $\psi(S_k(x)\cap S_{m-k}(y))$, and the constant value $\varphi(\psi(x_{k+1}))\geq 0$ at all points of $\psi(S_{k+1}(x)\cap S_{m-k-1}(y))$. Denoting $\gamma=\abs{\varphi(\psi(x_k))}>0$, one concludes: the closed convex hulls of $\psi(S_k(x)\cap S_{m-k}(y))$ and of $\psi(S_{k+1}(x)\cap S_{m-k-1}(y))$ are at least $\gamma>0$ apart. 

Our choice of $m$ assures that $k\geq 2$, and in particular the intersections $S_k(x)\cap S_{m-k}(y)$ and $S_{k+1}(x)\cap S_{m-k-1}(y)$ are infinite. In fact, they are both isometrically isomorphic to the Urysohn metric space $\Ur_{0,1,\ldots,k}$ of diameter $2k\geq 4$.

Let $N$ be given. Choose $2N+2$ points $a_i,b_i$, $i=0,1,2,\ldots,N$ on the sphere $S_m(z_0)$ so that $a_0=x$, $b_0=y$ and the distances between any two distinct points from among them is given by:
\[d(z,w)=\left\{\begin{array}{ll} 1,&\mbox{ if }z=a_i,~w=a_j,~i\neq j,\\
1,&\mbox{ if }z=b_i,~w=b_j,~i\neq j,\\
m,&\mbox{ if }z=a_i,~w=b_i, \\
m-1,&\mbox{ if }z=a_i,~w=b_j,~~i\neq j.
\end{array}\right.\]
For an arbitrary sequence $\e=(\e_i)_{i\leq N}\in\{0,1\}^N$, define the function $f_\e$ on $\{z_0\}\cup \{a_i\}_{i\leq N}\cup \{b_i\}_{i\leq N}$ by the conditions:
\[f_\e(z_0)=m,~~f_\e(a_i) = k+\e_i,~~f_\e(b_j)=m-k-\e_j.\]
Now one can verify, by considering 17 separate cases, that $f_\e$ is an (integer-valued) Kat\v etov function and so a distance function from some point $x^\ast\in S_m(z_0)$. 
It means that the intersection
\[T_\e=\cap_{i=0}^N S_{k+\e_i}(a_i)\cap\cap_{i=1}^N S_{m-k-\e_i}(b_i)\cap S_m(z_0)\]
is non-empty. 

Associate to every $\e\in\{0,1\}^N$ the closed convex hull $C_\e$ of $\psi(T_\e)$, which is a weakly compact subset of $F$. If $\e\leq\delta$, $\e\neq\delta$, then $C_\e$ and of $C_\delta$ are at a distance of at least $\gamma>0$ from each other, where the constant $\gamma$ was defined previously in this proof. Indeed, suppose $0\leq i\leq N$ be such that $\e_i\neq\delta_i$, and
let $h$ be an isometry of $\Ur_\Z$, preserving $z_0$ and taking $x\mapsto a_i$ and $y\mapsto b_i$.
The linear functional $\varphi\circ h$ has norm one and assumes constant values on $C_\e$ and on $C_\delta$, differing between themselves by at least $\gamma$.

Denote by $\mathscr T$ the prefix binary tree associated to the Hamming cube $\{0,1\}^N$, i.e., $\mathscr T$ consists of all prefix strings of elements $\e\in\{0,1\}^N$, with $\sigma\leq\tau$ if and only if $\sigma$ is a prefix of $\tau$. The system $C_\sigma$, $\sigma\in {\mathscr T}$ forms a binary tree under inclusion, and for every two nodes $\sigma,\tau$ at the same level the distance between $C_\sigma$ and $C_\tau$ is at least $\gamma$.  
By Corollary \ref{c:corol}, the space $F$ is not uniformly convex.
\qed
\end{pf}

\section{Appendix: alternative proof of Solecki-Vershik theorem \ref{th:sv}}

As brought to this author's attention by S. Solecki, the first version of his article \cite{solecki} contained a proof of theorem \ref{th:sv} along the same lines as outlined below. We include the proof just to make the paper reasonably self-contained. 

Let $X$ be a finite metric space. Denote by $P$ the set of all partial isometries $p$ of $X$ whose domain ${\mathrm{dom}}\, p$ is non-empty.
Let $F=F(P)$ be the free group on $P$. 
Every word $w\in F$ defines in a unique way a partial isometry of $X$ (possibly one with empty domain), under the convention that the empty word $e$ corresponds to the identity map of $X$. In this way, one obtains a {\em partial action} of $F(P)$ on $X$, that is, a map from $F(P)$ to the set of partial isometries of $X$ satisfying the properties that for all $x\in X$:
\begin{enumerate}
\item $e\cdot x=x$,
\item if $u\cdot x$ is defined, then $u^{-1}\cdot u \cdot x=x$,
\item if $u\cdot v\cdot x$ is defined, then $(uv)\cdot x = u\cdot v\cdot x$.
\end{enumerate}

Cf. \cite{exel,MS}.

A {\em globalization} of a pair consiting of a metric space $X$ and a partial action by a group $G$ on $X$ is a metric space $Y$ containing $X$ as a metric subspace and equipped with a global action of $G$ in such a way that for every $g\in G$ the partial isometry of $X$ defined by $g$ is a restriction of the corresponding global isometry of $Y$.
A {\em universal globalization} of $X$ is a globalization $Z$ with the property that the embedding of $X$ into any other globalization uniquely factors through the embedding $X\hookrightarrow Z$. 

Given a finite metric space $X$, we will consider it as equipped with a canonical partial action of the free group $F=F(P)$. Since every finite subspace is $g$-embedded into the Urysohn metric space, every isometric embedding $X\hookrightarrow \Ur$ determines (in a more than one way) a globalization of $X$. 

The universal globalization of a pair $(X,F(P))$ as above was constructed by Megrelishvili and Schr\"oder in \cite{MS}. 
We will denote this globalization by $U(X)$. 
The construction of a globalization of a partial action appears in many previous works and, as stressed by the anonymous referee, goes back at least to Mackey. A contribution of \cite{MS} is, in particular, establishing that the canonical mapping from $X$ to $U(X)$ is an isometry under a rather weak set of assumptions. 

At the set level, this globalization is the quotient set of $F\times X$ modulo the equivalence relation 
\begin{equation}
(uv,x)\sim (u, v\cdot x)\mbox{ whenever }v\cdot x\mbox{ is defined.}
\label{eq:MS}
\end{equation}
As in \cite{MS}, we will denote the equivalence class of a pair $(u,x)$ by $[u,x]$. The action of $F$ on $U(X)$ is defined by $g\cdot [u,x]=[gu,x]$. 

The universal globalization $U(X)$ admits an alternative description as a homogeneous factor-space of the group $F(P)$. Here is a repetition,
{\em mutatis mutandis,} of a construction presented in \cite{HL} on pp. 1987--1988. 
Choose a point $a_0\in X$ and denote by $H_0$ a subgroup of $F$ generated by the set
\begin{eqnarray*}
X_0&=&\{p^{-1}\cdot p^\prime\colon p,p^\prime\in P,~~p(a_0)=p^\prime(a_0)\}\cup \\
&& \{p_3^{-1}\cdot p_1\cdot p_2\colon p_1,p_2,p_3\in P,~~p_1\circ p_2(a_0)=p_3(a_0)\}.\end{eqnarray*}

Let now $H$ be any subgroup of $F$ satisfying $H_0<H$.

For every $a\in X$ there is a partial isomorphism $p\in P$ taking $a_0$ to $a$ (for instance, one with ${\mathrm{dom}}\,p=\{a_0\}$, ${\mathrm{im}}\,p=\{a\}$). Furthermore, if $p^\prime$ has the property $p^\prime(a_0)=a$, then the left cosets $pH$ and $p^\prime H$ coincide. Therefore, the map $\phi$ from $X$ to the homogeneous space $F/H$ given by the formula
\[\phi(a) = pH,\mbox{ where }p\in P\mbox{ and }p(a_0)=a,\]
is well-defined.  
If moreover $H\cap X_1=\emptyset$, where
\[X_1=\{p^{-1}\cdot p^\prime\colon p,p^\prime\in P,~~p(a_0)\neq p^\prime(a_0)\},\]
then $\phi$ is injective.
We will assume this condition to be satisfied, and will identity $X$ with its image under $\phi$ in $F/H$. 

Every $g\in F$ determines a left translation of the factor-space $F/H$, which we will denote $\tilde g$. It is easy to see that for every $p\in P$ and every $a\in{\mathrm{dom}}\,p$ one has $\tilde p(a)=p(a)$. Indeed, this condition means, in full, $p\phi(a)=\phi(p(a))$, or $pp_1H=p_2H$, where $p_1(a_0)=a$ and $p_2(a_0)=p(a)$. Since $p_2^{-1}pp_1(a_0)=a_0$, one has $p_2^{-1}pp_1\in H_0\subseteq H$, and the condition holds. 

In order to make $F/H$ into a metric space, we first turn it into an edge-coloured graph. Namely, we add an edge labeled with a real number $r$ between two elements $\alpha$ and $\alpha^\prime$ if and only if there are $a,a^\prime\in X$ and a $g\in F$ with 
\begin{equation}
d(a,a^\prime)=r,~~ga=\alpha,\mbox{ and }ga^\prime=\alpha^\prime.
\label{eq:prime}
\end{equation}
(Again, we identify $a$ with the corresponding coset $\phi(a)$, etc.)

Since for every word $w\in F(P)$ and each $p\in P$ the cosets $wH$ and $wpH$ are adjacent, with the corresponding edge carrying the weight $d(a_0,p(a_0))$, an inductive argument shows that the graph built on $F/H$ is connected.
Now we equip $F/H$ with the corresponding path distance, which is a left-invariant pseudometric.

There are a few potential problems that may arise here. Firstly, is the edge-labeling as above uniquely defined? Secondly, is the path distance a genuine metric (that is, the distance between two distinct points is non-zero)? Thirdly, will the restriction of this distance to $X$ coincide with the original metric on $X$? (Apriori, it is only bounded by $d_X$ from above.) 

In the case where $H=H_0$,
the answer to all three questions is positive, and it follows from the construction of Megrelishvili and Schr\"oder mentioned above. 

Indeed, consider the following formula for an arbitrary $g\in F$:
\[\psi(gH_0) = [g,a_0].\]
Since every element of $H_0$ stabilizes $a_0$, one has $[gh,a_0]=[g,h(a_0)]=[g,a_0]$, and the map $\psi$ from $F/H_0$ to $U(X)$ is well-defined. Clearly, $\psi$ is surjective and $F$-equivariant. Further, it is not difficult to verify that $H_0$ is precisely the stabilizer of the class $[e,a_0]\in U(X)$, and so the map $\psi$ is a bijection. Lifting the metric from $U(X)$ to $F$, one concludes that the path distance constructed above is a metric extending the distance $d_X$.

All that remains to be done, is to show that the same three conclusions hold for at least one subgroup $H<F$ of finite index containing $H_0$. Notice that the condition that the path distance, $d$, be a metric on $F/H$ is not essential: as long as the restriction of $d$ to $X$ coincides with the distance on $X$, one can replace $(F/H,d)$ with the associated metric space, that is, the quotient space under the relation $x\sim y\iff d(x,y)=0$. 

The two conditions to be verified are that the edge labeling be uniquely defined, and that 
the path distance on $F/H$ be an extension of the distance $d_X$ on $X$. Both can be reformulated in a unified way as follows. Let $a,b\in X$. Say that two finite sequences of pairs $(c_i,d_i)$ of elements of $X$ and of elements $x_i$ of $F$, $i=1,2,\ldots,n$, form a {\em bad configuration} for $(a,b)$, if the following conditions are met: 
\begin{enumerate}
\item $x_1\cdot c_1=a$,
\item $x_{i+1}c_i=x_id_i$ for $i=1,2,\ldots,n$,
\item $x_nd_n=b$, and
\item $\sum_{i=1}^n d_X(c_i,d_i)<d_X(a,b)$.
\end{enumerate} 

In particular, non-existence of bad configurations implies the uniqueness of labeling.

We want to stress again that here we identify $X$ with its image in $F/H$ under $\phi$. A more accurate rendering of the existence of a bad configuration in $F/H$ is therefore given by the following {\em ad hoc} concept.

\begin{definition}
\label{d:bad}
Let $X$ be a finite metric space and $H$ a subgroup of $F(P)$ containing $H_0$.
A {\em bad configuration} for $X$ modulo $H$ is a collection of
elements $p,q,p_i,q_i\in P$ and $x_i\in F(P)$,  $i=1,2,\ldots,n$, such that 
\begin{enumerate}
\item $x_1p_1\equiv p \,{\mathrm{mod}}\,H$,
\item $x_{i+1}p_i\equiv x_iq_i\,{\mathrm{mod}}\,H$ for $i=1,2,\ldots,n$,
\item $x_nq_n\equiv q\,{\mathrm{mod}}\,H$, and
\item $\sum_{i=1}^n d_X(p_i(a_0),q_i(a_o))<d_X(p(a_0),q(a_0))$.
\end{enumerate}
\end{definition}

What we want, is to avoid bad configurations for $X$ modulo $H$ by carefully choosing a finite index subgroup $H$ that in addition contains $H_0$. Remark that there are only finitely many theoretically possible bad configurations for any given finite metric space $X$, so if we can learn to choose a subgroup $H$ of finite index containing $H_0$ and avoiding a given bad configuration, we are done, as the intersection of finitely many subgroups of finite index has finite index by Poincar\'e's theorem.

Recall that a group $G$ is {\em residually finite} if homomorphisms from $G$ to finite groups separate points. For example, free groups are residually finite. Equivalently, residual finiteness means that for every finite subset $A\subseteq F\setminus\{e\}$, where $F$ denotes a non-abelian free group, there is a normal subgroup $N\triangleleft F$ of finite index disjoint from $A$. In fact, the collection of subgroups of free groups of finite index is much richer than that, as shows the following surprising result.

\begin{theorem}[M. Hall \cite{hall}] Every finitely generated subgroup $H$ of a non-abelian free group $F$ is the intersection of subgroups $N<F$ of finite index. \qed
\end{theorem}

This result can be further strengthened.

\begin{theorem}[Ribes and Zalesski\u\i\ \cite{rz}] Let $H_1,H_2,\ldots,H_n$ be finitely generated subgroups of a non-abelian free group $F$. Then for every $g\in F\setminus H_1H_2\ldots H_n$ there are finite index subgroups $K_1,K_2,\ldots, K_n<F$ such that $H_i<K_i$, $i=1,2,\ldots,n$, and $g\notin K_1K_2\ldots K_n$. \qed
\end{theorem}

A further refinement of the above theorems forms the core result of the above mentioned paper by Herwig and Lascar. To state it, we need to recall their terminology from \cite{HL}. A {\em left-system} is a finite set of equations of the form
\[x\equiv_i y\cdot g\mbox{ or }x\equiv_i g,\]
where $i=1,2,\ldots,n$, $x,y$ belong to a finite set of unknowns, $X$, and $g$ are elements of a free non-abelian group $F$. If ${\mathcal H}=(H_1,H_2,\ldots, H_n)$ is a sequence of subgroups of $F$, a solution of a left-system as above modulo $\mathcal H$ is a family $g_x$, $x\in X$ of elements of $F$ such that for every equation of the form $x\equiv_i y\cdot g$ one has
\[g_x\equiv g_y\cdot g\,{\mathrm{mod}}\, H_i,\]
and for every equation of the form $x\equiv_i g$ one has
\[g_x\equiv g\,{\mathrm{mod}}\, H_i.\]

\begin{theorem}[Herwig and Lascar \cite{HL}] Let $n\in\N$, let ${\mathcal H}=(H_1,H_2,\ldots, H_n)$ be a sequence of subgroups of a free non-abelian group $F$, and let $(E)$ be a left-system of equations in $F$. Assume that $(E)$ has no solutions in $F$ modulo $\mathcal H$. Then there exist finite index subgroups $K_1,K_2,\ldots,K_n$ of $F$ such that $H_i<K_i$, $i=1,2,\ldots,n$, and the left-system $(E)$ has no solutions in $F$ modulo ${\mathcal K}=(K_1,K_2,\ldots,K_n)$. \qed
\label{th:hl}
\end{theorem}

Notice that the theorem of Ribes and Zalesski\u\i\ is obtained from the above result if one considers the left-system of the form
\[\left\{\begin{array}{lll}
x_n&\equiv_n &g\\
x_{n-1}&\equiv_{n-1}&x_n\\
\vdots & \vdots & \vdots  \\
x_2&\equiv_2 & x_3 \\
x_1&\equiv_1& e.
\end{array}\right.\]
This system has no solutions in $F$ modulo $\mathcal H$ if and only if $g\notin H_1H_2\ldots H_n$. 

Now it remains to notice that the existence of a bad configuration for a finite metric space $X$ modulo a subgroup $H$ (Definition \ref{d:bad}) is equivalent to the existence of a solution to a left-system of equations:

\[\left\{\begin{array}{lll}
x_1&\equiv& pp_1^{-1} \\
x_2&\equiv& x_1\cdot q_1p_1^{-1} \\
x_3&\equiv & x_2\cdot q_2p_2^{-1} \\
\vdots & &\vdots \\
x_n&\equiv & qq_n^{-1}.
\end{array}\right.\]
Since there is no solution of this left-system modulo $H_0$, there is a subgroup $H$ of finite index containing $H_0$ and such that there are no solution modulo $H$ either. This finishes the proof of Theorem \ref{th:sv} for the case of the Urysohn space.

Note that the above argument in effect leads to the following slight technical generalization.
As usual, we say that a subset $S$ of a totally ordered set $X$ is {\em convex} if for every $x,y,z\in X$ the conditions $x\leq y\leq z$ and $x,z\in S$ imply $y\in S$. 

\begin{theorem}[General form of Hrushovski--Solecki--Vershik theorem]
\label{th:hsv}
Let $S$ be a convex subset of an additive subsemigroup $T$ of $\R$, containing zero. Then the universal Urysohn metric space $\Ur_S$ whose distance takes values in $S$ has the Hrushovski--Solecki--Vershik property. 
\end{theorem}

Remark that under our assumptions on $S$, the space $\Ur_S$ always exists and is unique, according to Theorem 1.4 and Example 1.5.3 in \cite{DLPS}. Now in the case $S=\{0,1,2\}$ (a convex subset of the semigroup $\N$ of natural numbers) one recovers Hrushovski's theorem \ref{th:h}, while the case $S=\R$ gives the Solecki--Vershik theorem \ref{th:sv}.

\section{Open questions}

\numero
The most interesting open question at the moment of writing this article seems to be whether or not every finitely generated group admits a coarse embedding into a superreflexive Banach space.

\numero It remains unknown whether a direct sum of graphs forming an expander family can ever admit a coarse embedding into a superreflexive Banach space (cf. open problem $\sharp$ 9, submitted by Piotr Nowak on the list \cite{open}). 

\numero Is there a proof of the Solecki--Vershik property along the lines of a simple combinatorial proof of the Hrushovski theorem given in \cite{HL}, Sec. 4.1? 

\numero
Is the following ``coarse analogue'' of Holmes' theorem \cite{holmes} true? Suppose the Urysohn space $\Ur$ is coarsely embedded into a normed space $E$ in such a way that the image of $\Ur$ spans $E$. Then $E$ is coarsely equivalent to the Lipschitz-free Banach space over $\Ur$. (Cf. also \cite{P06}, p. 112, as well as \cite{melleray}.)

\numero
Does there exist an analogue of the universal Urysohn metric space in the coarse category among spaces of bounded geometry and exponential growth? Cf. related constructions in \cite{dgly}.

\numero A {\em regular embedding} of a (simple, non-oriented) graph $\Gamma$ into a metric space $X$ is a map from the set of vertices of $\Gamma$ to $X$ such that the distance between images of two adjacent vertices is always $\alpha$, and between non-adjacent ones is always $\beta$, where $\alpha<\beta$. It is well-known that many finite graphs do not admit regular embeddings into the Hilbert space. (Cf. \cite{maehara}.) 
Does the infinite random graph $R$ admit a regular embedding into a reflexive Banach space? 

\section*{Acknowledgements}

The author is grateful to Christian Rosendal for spotting a serious error in the first version of the article, and to the anonymous referee for a number of insightful comments which have led to a considerable reworking and shortening of the article.
The research was supported by the University of Ottawa internal research
grants and by the 2003-07 NSERC discovery grant 
``High-dimensional geometry and topological transformation groups.''


\begin{thebibliography}{100}

\bibitem{BaL} F. Baudier and G. Lancien, {\em Embeddings of locally finite metric spaces into Banach spaces,} arXiv:math/0702266v1, to appear in Proc. Amer. Math. Soc.

\bibitem{BL} Y. Benyamini and J. Lindenstrauss,
\textit{Geometric Nonlinear Functional Analysis,} Vol. 1,
Colloquium Publications \textbf{48}, American Mathematical Society,
Providence, RI, 2000.

\bibitem{Blum} L.M. Blumenthal, 
{\it Theory and Applications of Distance Geometry,}
Chelsea Publ. Co., Bronx, NY, 1970.

% \bibitem{BG} N. Brown and E. Guentner, \textit{Uniform embeddings of bounded geometry spaces into reflexive Banach space,} Proc. Amer. Math. Soc. \textbf{133} (2005), 2045--2050.
 
\bibitem{cameron} P. Cameron, \textit{The random graph,} in:
The Mathematics of Paul Erdos (J. Ne\v setril, R. L. Graham, eds.),
Springer, 1996, pp. 331-351.

\bibitem{day}  M.M. Day, \textit{Some more uniformly convex spaces,}
Bull. Amer. Math. Soc. \textbf{47} (1941), 504--507. 

\bibitem{DLPS}
C. Delhomme, C. Laflamme, M. Pouzet, N. Sauer,
\textit{Divisibility of countable metric spaces}, European Journal of Combinatorics \textbf{28} (2007), 1746-1769.

\bibitem{dgly}
A.N. Dranishnikov, G. Gong, V. Lafforgue, and G. Yu, 
\textit{Uniform embeddings into {H}ilbert space and a question of {G}romov}, Canad. Math. Bull. \textbf{45} (2002), 60--70.
	 
\bibitem{enflo} P. Enflo, \textit{On a problem of Smirnov,} Ark. Mat. \textbf{8} (1969), 107--109. 

\bibitem{exel} R. Exel, {\em Partial actions of groups and actions of inverse semigroups,} Proc. Amer. Math. Soc. \textbf{126} (1998), 3481--3494.

\bibitem{fhhmpz}
M. Fabian, P. Habala, P. H\'ajek, V. Montesinos Santaluc'a, J. Pelant, V. Zizler,
\textit{Functional Analysis and Infinite-dimensional Geometry}, 
CMS Books in Mathematics, \textbf{8}, Springer-Verlag, New York, 2001. 

\bibitem{gromov91}
M. Gromov, {\em
Asymptotic invariants of infinite groups,} in: Geometric group theory, Vol. 2 (Sussex, 1991), 1--295, 
London Math. Soc. Lecture Note Ser., \textbf{182}, Cambridge Univ. Press, Cambridge, 1993.

\bibitem{Gr} M. Gromov, \textit{Metric Structures for Riemannian and
Non-Riemannian Spaces,} Progress in Mathematics \textbf{152}, Birkhauser
Verlag, 1999.

\bibitem{gromov00} M. Gromov, \textit{Spaces and questions,}
GAFA 2000 (Tel Aviv, 1999), 
Geom. Funct. Anal. 2000, Special Volume, Part I, 118--161. 

\bibitem{hall} M. Hall, {\em A topology for free groups and related groups,} Ann. of Math. \textbf{52} (1950), 127--139.

\bibitem{HI} C.W. Henson and J. Iovino,
\textit{Ultraproducts in analysis}, in: {Analysis and Logic (Mons, 1997)},
{London Math. Soc. Lecture Note Ser.}, \textbf{262}, {1--110},
{Cambridge Univ. Press}, {Cambridge}, {2002}.

\bibitem{HL}
B. Herwig and D. Lascar, {\em
Extending partial automorphisms and the profinite topology on free groups,}
Trans. Amer. Math. Soc. \textbf{352} (2000), 1985--2021. 

\bibitem{holmes} M.R. Holmes, \textit{The universal separable metric space of Urysohn and isometric embeddings thereof in Banach spaces,} Fund. Math. \textbf{140} (1992), 199--223. 

\bibitem{hrushovski}
E. Hrushovski, {\it Extending partial isomorphisms of graphs,}
Combinatorica {\bf 12} (1992), 411--416. 

\bibitem{jr}
W.B. Johnson and N.L. Randrianarivony, 
{\em $l\sb p (p>2)$ does not coarsely embed into a Hilbert space,}
Proc. Amer. Math. Soc. \textbf{134} (2006), 1045--1050.

\bibitem{kalton} N.J. Kalton, {\em Coarse and uniform embeddings into reflexive spaces,} to appear in Quart. J. Math. Oxford.

\bibitem{ky}
G. Kasparov and G. Yu, 
\textit{The coarse geometric Novikov conjecture and uniform convexity,}
arXiv math.OA/0507599,  to appear in Adv. Math.

\bibitem{kechris}
A.S. Kechris, \textit{Global Aspects of Ergodic Group Actions and Equivalence Relations,} preprint, November 10, 2006, 197 pp.

\bibitem{lafforgue} V. Lafforgue, {\em Expanders with no uniform embedding in a uniformly convex Banach space,} abstract of talk at the conference on Geometric Linearization of Graphs and Groups, Lausanne, January 22-26, 2007, 
\newline
{\tt http://bernoulli.epfl.ch/graphs/listeConfJanvier.php} (accessed Jan. 3, 2007).

\bibitem{proceedings} A. Leiderman, V. Pestov, M. Rubin, S. Solecki, V.V. Uspenskij (eds.), 
\textit{Proceedings of the Workshop on the Universal Urysohn Metric Space} (Ben-Gurion University of the Negev, 21-24 May 2006),
to appear as a Special Volume of Topology and its Applications. 

\bibitem{lw} J. Li and Q. Wang, 
\textit{Remarks on coarse embeddings of metric spaces into uniformly convex Banach spaces,}
J. Math. Anal. Appl. \textbf{320} (2006), 892--901. 

\bibitem{lanvt} J. L\'opez-Abad and L. Nguyen Van Th\'e, \textit{The oscillation stability problem for the Urysohn sphere: a combinatorial approach,} preprint, November 2006, 21 pp., to appear in \cite{proceedings}.

\bibitem{maehara}
H. Maehara, \textit{Regular embeddings of a graph,}
Pacific J. Math. \textbf{107} (1983), 393--402. 

\bibitem{megrelishvili01} M.G. Megrelishvili, \textit{Every semitopological semigroup compactification of the group $H\sb +[0,1]$ is trivial,} Semigroup Forum \textbf{63} (2001), 357--370.

\bibitem{megrelishvili07} M. Megrelishvili, {\em Topological transformation groups: selected topics,} in: Open Problems in Topology (Second Ed., Elliott Pearl, ed.), to be published by Elsevier Science, 2007.

\bibitem{MS} M. Megrelishvili and L. Schr\"oder, {\em Globalization of confluent partial actions on topological and metric spaces,} Topology Appl. \textbf{145} (2004), 119--145.

\bibitem{melleray}
J. Melleray, \textit{Some geomeric and dynamical properties of the Urysohn space,} 2006 preprint, 
a contribution to the volume \cite{proceedings}.

\bibitem{nowak05}
P.W. Nowak, 
\textit{Coarse embeddings of metric spaces into Banach spaces,}
Proc. Amer. Math. Soc. \textbf{133} (2005), 2589--2596.

\bibitem{nowak06}
P.W. Nowak, \textit{On coarse embeddability into $l\sb p$-spaces and a conjecture of Dranishnikov,} Fund. Math. \textbf{189} (2006), 111--116. 

\bibitem{open} \textit{Open problems raised at the Concentration Week on Metric Geometry and Geometric Embeddings of Discrete Metric Spaces,} URL: \hfill\newline
{\tt http://www.math.tamu.edu/research/workshops/linanalysis/problems.html}
(accessed on 27 December 2006). 

% \bibitem{P02} V. Pestov, \textit{Ramsey--Milman phenomenon, 
% Urysohn metric spaces,
% and extremely amenable groups,} Israel Journal of Mathematics
% \textbf{127} (2002), 317-358. \textit{Corrigendum,}  
% ibid., \textbf{145} (2005), 375-379.

\bibitem{P06}
V. Pestov, \textit{Dynamics of Infinite-dimensional Groups: the Ramsey--Dvoretzky--Milman phenomenon,} 
University Lecture Series {\bf 40}, American Mathematical Society, Providence, RI,  2006.

\bibitem{Rado} R. Rado, \textit{Universal graphs and universal functions,}
Acta Arithm. \textbf{9} (1954), 331--340.

\bibitem{rz} L. Ribes and P.A. Zalesski\u\i, {\em On the profinite topology of the free group,} Bull. London Math. Soc. \textbf{25} (1993), 37--43.

\bibitem{roe} J. Roe, \textit{Lectures on Coarse Geometry,} University Lecture Series \textbf{31}, American Mathematical Society, Providence, RI, 2003. 

\bibitem{solecki} S. Solecki, \textit{Extending partial isometries,} Israel J. Math. \textbf{150} (2005), 315--332.

\bibitem{Usp90}
V.V. Uspenski\u\i,  {\it On the group of isometries of the
Urysohn universal metric space,} Comment. Math. Univ. Carolinae
{\bf 31} (1990), 181-182.

\bibitem{Usp98} V.V. Uspenskij, {\it On subgroups of minimal
topological groups,} 1998 preprint, later prepublished at
arXiv:math.GN/0004119, a contribution to the volume \cite{proceedings}.

\bibitem{Ver98} A.M. Vershik, \textit{The universal Urysohn space, 
Gromov's metric triples, and random metrics on the series of positive 
numbers,} Russian Math. Surveys \textbf{53} (1998), 921--928.

\bibitem{Ver04} A.M. Vershik, \textit{The universal and random metric
spaces,} Russian Math. Surveys \textbf{356} (2004), 65--104.

\bibitem{Ver05b} A.M. Vershik, \textit{Extensions of the partial isometries
of the metric spaces and finite approximation of the group of isometries of
Urysohn space,} preprint, 2005.

\end{thebibliography}
\end{document}